\documentclass{article}

\newcommand{\be}{\begin{equation}}
\newcommand{\ee}{\end{equation}}

\newcommand{\D}{\displaystyle}

\newcommand{\beq}{\begin{eqnarray}}
\newcommand{\eeq}{\end{eqnarray}}
\newcommand{\nbeq}{\begin{eqnarray*}}
\newcommand{\neeq}{\end{eqnarray*}}

\begin{document}

 \title{\bf Characterizations of exponential distribution via conditional expectations 
 of record values}
\author{George P. Yanev \\
Department of Mathematics \\
University of Texas - Pan American \\
Edinburg, Texas 78539 \\
e-mail: yanevgp@utpa.edu}
\date{\empty}

\maketitle

\begin{abstract}
We prove that the exponential distribution is the only one which satisfies a regression identity. This identity involves
conditional expectation of the sample mean of record values given two record values outside of the sample.
\end{abstract}

\vspace{0.3cm}\noindent {\it MSC Classification:} 62G30; 62E10.

\vspace{0.3cm} \noindent {\it Keywords:} characterizations, exponential distribution, record values.

\section{Discussion of results}
\label{intro}

 Before formulating and discussing
 the obtained results, let us introduce some notation regarding record values. Let $X_1, X_2,
\ldots$ be independent copies of a random variable $X$ with
absolutely continuous cumulative distribution function (c.d.f.) $F(x)$. (Upper) Record value is an observation
in a discrete time series, which exceeds all previous
observations, i.e., $X_j$ is a (upper) record value if $X_j>X_i$ for all
$i<j$. Define the sequence $\{T_n, n\ge 1\}$ of record times by $T_1=1$ and $
T_n=\min \{j:X_j>X_{T_{n-1}}, j>T_{n-1}\}$ for $n>1$. Then the
corresponding record values are $R_n=X_{T_n}$ for  $n=1,2,\ldots$ (see Nevzorov (2001)).

Let $F(x)$ be the exponential c.d.f. given by
\begin{equation}
\label{exp_type}
F(x)=1-e^{\D -cx}, \qquad (x\ge 0),
\end{equation}
where $c>0$. It is well-known (cf. Arnold et al. (1998), Section 4.2.2), that a single regression function (conditional expectation) of the form $E(R_{n+1}|R_n)$ is enough to characterize a continuous c.d.f. Several authors have studied characterizations of exponential
distributions in terms of the regression of one record value with
two other record values as covariates, i.e., regression functions given for $1\le k\le n-1$ and
$r\ge 1$ by
\be \label{general} E[g(R_n)|R_{n-k}=u, R_{n+r}=v ]=\varphi(u,v), \qquad (0\le u<v),\ee
where the function $g$ satisfies certain regularity
conditions.
Let $\bar{g}(u,v)$ denote the average value of the function $g(x)$ over the interval $[u,v]$, i.e.,
\[
\bar{g}(u,v)=\frac{1}{v-u}\int_u^v g(t)dt.
\]
Bairamov et al. (2005) study the particular case of (\ref{general}) when both
covariates are adjacent (only one spacing away) to $R_n$. They prove, under some
regularity conditions on $F$ and $g$ and using different notation, that  $F$ is exponential  if and only if
\[
E\left[g(R_n){\Big |}R_{n-1}=u, R_{n+1}=v \right]=\bar{g}(u,v), \qquad (0\le u<v).\]
On the other hand, Yanev et al. (2008) prove the following necessary condition for a continuous c.d.f. to be exponential. If $F$ is exponential then for $1\le k\le n-1$ and
$r\ge 1$
\be \label{exp_property}
\hspace{-0.5cm}E\left[\frac{\D g^{(k-1+r-1)}(R_n)}{\D k+r-1}{\Big |}R_{n-k}=u,
R_{n+r}=v\right]
=
{k-1+r-1 \choose k-1}\frac{\partial^{k-1+r-1}}{\partial u^{r-1}\partial
v^{k-1}}\ \bar{g}(u,v), \nonumber \ee
where $0\le u<v$. Let us demonstrate (\ref{exp_property}) using two particular choices of
the function $g(x)$. If $g(x)=x^{k+r-1}/(k+r-1)!$, then (\ref{exp_property}) becomes
\be \label{YAB08c}
E[R_n|R_{n-k}=u, R_{n+r}=v]=\frac{ru+kv}{k+r}, \qquad (0\le u<v).
\ee
Note that the right-hand side of (\ref{YAB08c}) is
the weighted mean of the covariates' values. The weight of one covariate
is proportional to the number of spacings $R_n$ is away from the other covariate.
As a result, the closer $R_n$ is to one of the covariates, the more affected its value is from this covariate.
Another choice of the function $g(x)$ could be
$g(x)=(-1)^{k+r-2}x^{-2}/(k+r-1)!$ when (\ref{exp_property}) is equivalent to
\[
E\left[\frac{1}{R^{k+r}_n}{\Big |}R_{n-k}=u, R_{n+r}=v\right]=\frac{1}{u^rv^k}, \ \  (0\le u<v).
\]
Whether (\ref{exp_property}) is not only a necessary, but also a sufficient condition for an absolutely continuous $F$ to be exponential for any $1\le k\le n-1$ and $r\ge 1$ is an open question. Yanev et al. (2008) obtain an affirmative answer when $r=1$, i.e., one
of the two covariates is adjacent to $R_n$. They show,  under some regularity assumptions on $F$ and $g$, that $F(x)$ is
exponential if and only if
\be \label{YAB08} E\left[\frac{g^{(k-1)}(R_n)}{k}{\Big |}R_{n-k}=u,
R_{n+1}=v\right]= \frac{\partial^{k-1}}{\partial
v^{k-1}}\ \bar{g}(u,v),  \qquad (0\le u<v). \ee

A different direction for extending the characterization results above is to consider conditional expectations of sums of record values.
Akhundov and Nevzorov (2007) prove that an absolutely continuous $F$ is exponential if and only if
\beq \label{mean_x}
E\left[\frac{R_2+R_{3}+\ldots+R_{n}}{n-1}| R_{1}=u, R_{n+1}  =  v\right] & =  & \frac{u+v}{2}.
\eeq
Our first result is one generalization of (\ref{mean_x}).

\vspace{0.5cm}{\bf Theorem 1}\ {\it Let $n\ge 2$ be an integer. Suppose $F(x)$ is absolutely continuous, $g(x)$ is continuous
in $[0 , \infty)$, and $g(v)\ne \bar{g}(0,v)$ for every $v>0$.
 Then (\ref{exp_type}) holds
 if and only if for $0\le u<v<\infty$
\be\label{mean}
E\left[\frac{g(R_2)+g(R_3)+\ldots+g(R_n)}{n-1}| R_{1}=u, R_{n+1}  =  v\right] =   \bar{g}(u,v).
\ee
}
Next, we extend Theorem 1 to the case when the second covariate is two spacings away from $R_n$. The general case of non-adjacent covariates remains an open problem. Denote, for a continuous
 function $g(x)$ and  $n=1,2,\ldots$
\be
 \label{I_n}
I_n(u,v)=\int_0^1 g(u+(v-u)z)z^{n-1}dz,\qquad (0\le u<v).
\ee
Note that $I_1(u,v)=\bar{g}(u,v)$.
Denote
$
H(x)=-\ln(1-F(x))$ for $x\ge 0$ and $h(x)=H'(x)$,
i.e., $H(x)$ is the cumulative hazard function of $X$ and $h(x)$ is its hazard (failure) rate.

\vspace{0.5cm}{\bf Theorem 2}\ {\it
Suppose $F(x)$ satisfies the following conditions.

(i) $F(x)$ is absolutely continuous in $[0, \infty)$ and $F''(x)$ is continuous in $(0,\infty)$;

(ii) $h(x)$ is nowhere constant in a small interval $(0, \varepsilon)$ for $\varepsilon>0$.

Assume that $g(x)$ is continuous in $[0, \infty)$, $g'(x)$ is continuous
in $(0 , \infty)$, and $g'(0+)\ne~0$. Let $n\ge 2$ be an integer.
Then (\ref{exp_type}) holds
if and only if for $0\le u<v<\infty$
\be \label{mean2}
\hspace{-0.3cm}E\left[\frac{g(R_2)+g(R_3)+\ldots+g(R_n)}{n-1}| R_{1}=u, R_{n+2}  =  v\right]
=\bar{g}(u,v) -\frac{nI_n(u,v)-\bar{g}(u,v)}{n-1}.
\ee
}

Setting in Theorem 2, $g(x)=x$ and thus $
nI_n(u,v)= (u+nv)/(n+1)$, and taking into account (\ref{mean2})
we obtain the following extension of (\ref{mean_x}).

\vspace{0.5cm}{\bf Corollary 1}\ {\it  Suppose $F(x)$ satisfies the following conditions.

(i) $F(x)$ is absolutely continuous in $[0, \infty)$ and $F''(x)$ is continuous in $(0,\infty)$;

(ii) $h(x)$ is nowhere constant in a small interval $(0, \varepsilon)$ for $\varepsilon>0$.
 Let $n\ge 2$ be an integer. Then (\ref{exp_type}) holds
 if and only if for $0\le u<v<\infty$
\beq \label{mean_x2}
E\left[\frac{R_2+R_{3}+\ldots+R_n}{n-1}| R_{1}=u, R_{n+2}  =  v\right]
 & = &   \frac{u+v}{2}-\frac{v-u}{2(n+1)}\\
 &  = &   \frac{(n+2)u+nv}{n+2+n}. \nonumber
\eeq
}
Notice that (\ref{YAB08c}) implies that a necessary condition for $F$ to be exponential is
\be \label{note}
E\left[\frac{R_2+R_{3}+\ldots+R_n}{n-1}| R_{1}=u, R_{n+m}  =  v\right]=\frac{(n+2m-2)u+nv}{n+2m-2+n},\qquad m\ge 1.
\ee
The results (\ref{mean_x}) and Corollary 1 prove the sufficiency of (\ref{note}) for $m=1$ and $m=2$.

\section{Proofs }

{\bf Lemma}\ Let $r$ and $n\ge 2$ be positive integers. If (\ref{exp_type}) holds, then
\beq \label{lemma}
\lefteqn{\hspace{-1.5cm}E\left[\sum_{j=2}^{n} g(R_j)\ |\ R_{1}=u, R_{n+r}=v\right]}\\
& & =(N+1)\int_0^1 g(u+(v-u)z) \left[ 1-\sum_{j=n-1}^{N}{N \choose j} z^j (1-z)^{N-j}\right]dz, \nonumber
\eeq
where $N=n+r-3$. As usual we assume $\sum_j^i (\cdot) =0 $ when $i<j$.

\vspace{0.3cm} {\bf Proof.}\
It is well-known (e.g., Nevzorov (2001), p.11 and p.69) that if (\ref{exp_type}) holds, then
$
R_n~=~\sum_{j=1}^n \xi_j,
$
where $\{\xi_j\}$ are i.i.d. exponential with (\ref{exp_type}). Therefore, for $2\le j\le n$
\nbeq
\lefteqn{E\left[ R_j| R_{1}=u, R_{n+r}=v\right]}\\
    & = &
E\left[ \xi_1+\xi_2+\ldots + \xi_j| \xi_{1}=u, \xi_1+\xi_2+\ldots +\xi_{n+r}=v\right] \\
    & = &
 E\left[ u+\xi_{2}+\ldots + \xi_j| \xi_{2}+\ldots +\xi_{n+r}=v-u\right]\\
    & = &
  E\left[ u+Y_{j}| Z_{n+r}=v-u\right],\qquad \mbox{say}.
\neeq
We obtain the following density functions related to $Y_j$ and $Z_{n+r}$ for $2\le j\le n$
\nbeq
f_{Y_j,Z_{n+r}}(y,z)& = & \frac{c^{n+r-1}y^{j-2}(z-y)^{n+r-j-1}}{(j-2)!(n+r-j-1)!}e^{-cz}, \quad
f_{Z_{n+r}}(z) =
\frac{c^{n+r-1}}{(N+1)!}z^{N+1}e^{-cz},\\
& & \\
f_{Y_j|Z_{n+r}}(y|z) 
    & = &
    \frac{(N+1)!}{(j-2)!(n+r-j-1)!}\left(\frac{y}{z}\right)^{j-2}\left(1-\frac{y}{z}\right)^{n+r-j-1}\frac{1}{z},
\neeq
where $N=n+r-3$. Now we are in a position to prove (\ref{mean}). Indeed,
\nbeq
\lefteqn{E\left[ \sum_{j=2}^{n} g(R_j)| R_{1}=u, R_{n+r}=v\right]
    =
    E\left[ \sum_{j=2}^{n} g(u+Y_j)| Z_{n+r}=v-u\right]}\\
    & = &
     \int_0^{v-u}\!\!\!\!g(u+t)\sum_{j=2}^{n} \frac{(N+1)!}{(j-2)!(n+r-j-1)!}\left(\frac{t}{v-u}\right)^{j-2}\!\!\!\!
    \left(1-\frac{t}{v-u}\right)^{n+r-j-1}\!\!\!\!\frac{dt}{v-u}\\
    & = &
     (N+1)\int_0^1 g(u+(v-u)z)\sum_{j=2}^{n} { N\choose j-2}z^{j-2}
    (1-z)^{n+r-j-1}dz\\
    & = &
    (N+1)\int_0^1g(u+(v-u)z)\sum_{i=0}^{n-2} { N\choose i+k-1}z^{i+k-1}(1-z)^{N-(i+k-1)}dz,
\neeq
which is equivalent to (\ref{lemma}).

{\bf \subsection{ Proof of Theorem 1 }}
Recall the formula for the conditional density $f_{j|1,n+r}(t|u,v)$, say, of $R_j$ given $R_{1}=u$
 and $R_{n+r}=v$, where $r\ge 1$ and $2\le j\le n$ (e.g., Ahsanullah (2004), p.6). It is given for $u<t<v$ and $2\le j\le n$ by
\beq \label{cond_density}
\lefteqn{f_{j|1,n+r}(t|u,v)}\\
& = &  \frac{(n+r-2)!}{(j-2)!(n+r-1-j)!}\left[\frac{H(t)-H(u)}{H(v)-H(u)}\right]^{j-2}
\left[\frac{H(v)-H(t)}{H(v)-H(u)}\right]^{n+r-1-j}\hspace{-0.4cm}\frac{h(t)}{H(v)-H(u)}\nonumber .
\eeq

\underline{Sufficiency}. Assuming (\ref{mean}), we will prove that $F(x)$ is exponential (\ref{exp_type}).
Using (\ref{cond_density}) with $r=1$ and the binomial formula, we obtain
\be \label{new_cond_den}
\sum_{j=2}^n f_{j|1,n+1}(t|u,v)=\frac{(n-1)h(t)}{H(v)-H(u)}.
\ee
Taking into account (\ref{new_cond_den}), for the left-hand side of (\ref{mean}) we have
\nbeq
E\left[ \sum_{j=2}^{n} g(R_j)| R_{1}=u, R_{n+1}=v\right] & = &
    \int_u^v \sum_{j=2}^{n} g(t)f_{j|1,n+1}(t|u,v)\ dt\\
        & = &
    \frac{n-1}{H(v)-H(u)} \int_{u}^v g(t)d H(t)
\neeq
and thus (\ref{mean}) is equivalent to
\[
\frac{1}{H(v)-H(u)} \int_{u}^v g(t)d H(t)=\frac{1}{v-u}\int_u^v g(t)dt.
\]
Letting $u\to 0+$, and hence $H(u)\to 0$, we obtain
\[
\int_{0}^v g(t)d H(t)=\frac{H(v)}{v}\int_{0}^v g(t)dt.
\]
Differentiating both sides with respect to $v$, after rearranging terms, we find
\[
[g(v)-\bar{g}(0,v)] h(v)=\frac{H(v)}{v}[g(v)-\bar{g}(0,v)].
\]
Dividing both sides by $g(v)-\bar{g}(0,v) \ne 0$,
we obtain the equation
\[
\frac{h(v)}{H(v)}=\frac{1}{v},
\]
which only solution is the exponential $F(x)$ given by (\ref{exp_type}). Indeed, integrating both sides with respect to $v$, we
obtain
\[
\ln H(x)=\ln v + \ln c    \qquad (c>0)
\]
and thus $ H(x)=-\ln (1-F(x))=cx$ and (\ref{exp_type}) follows.

\underline{Necessity}. If $F(x)$ is exponential, then (\ref{mean}) follows from the lemma setting $r=1$.

{\bf \subsection{ Proof of Theorem 2}}

\underline{Sufficiency}.
Recall the notation $f_{j|1,n+r}(t|u,v)$ introduced before (\ref{cond_density}). Now, (\ref{cond_density}) with $r=2$ and the binomial formula imply
\beq \label{new_identity}
\sum_{j=2}^n f_{j|1,n+2}(t|u,v) & = & \frac{nh(t)}{(H(v)-H(u))^{n}}\sum_{j=2}^n {n-1 \choose j-2}(H(t)-H(u))^{j-2}(H(v)-H(t))^{n-j+1}\nonumber \\
    & = &
    \frac{nh(t)}{(H(v)-H(u))^{n}}\sum_{i=0}^{n-2} {n-1 \choose i}(H(t)-H(u))^{i}(H(v)-H(t))^{n-1-i} \nonumber \\
    & = &
    \frac{nh(t)}{H(v)-H(u)}\left[ 1-\left(\frac{H(t)-H(u)}{H(v)-H(u)}\right)^{n-1}\right].
    \eeq
It follows from (\ref{new_identity}) that
\beq \label{mean22}
\lefteqn{E\left[ \sum_{j=2}^{n} g(R_j)| R_{1}=u, R_{n+2}  =  v\right]
     =
\int_u^v g(t)\sum_{j=2}^n f_{j|1,n+2}(t|u,v)dt} \\
& & =
\frac{n}{H(v)-H(u)}\int_u^v g(t)\left[ 1-\left(\frac{H(t)-H(u)}{H(v)-H(u)}\right)^{n-1}\right]dH(t). \nonumber
\eeq
It follows from (\ref{mean2}) and (\ref{mean22}), letting $u\to 0$ (and hence $H(u)\to 0$) and multiplying by $H^n(v)(n-1)/n$ that
\[
\int_{0}^v g(t)\left[H^{n-1}(v)-H^{n-1}(t)\right]dH(t)=H^n(v)[\bar{g}(v)-I_n(v)],
\]
where $\bar{g}(v)=\bar{g}(0,v)$ and $I_n(v)=I_n(0,v)$.
Differentiating with respect to $v$ and dividing by $H^{n-2}(v)h(v)$ we obtain
\[
(n-1)\int_{0}^vg(t)dH(t) = n H(v)[\bar{g}(v)-I_n(v)]+\frac{H^2(v)}{h(v)}[\bar{g}'(v)-I'_n(v)].
\]
Differentiating with respect to $v$ again and dividing by $H(v)$ we have
\beq \label{middle_eqn}
\lefteqn{\hspace{-1cm}[(n-1)g(v)-n\bar{g}(v)+nI_n(v)]\frac{h(v)}{H(v)}} \\
& & =\left[n+2-\frac{h'(v)H(v)}{h^2(v)}\right][\bar{g}'(v)-I'_n(v)]
+[\bar{g}''(v)-I''_n(v)]\frac{H(v)}{h(v)}. \nonumber
\eeq
Observing that
\be \label{identities}
nI_n(v)=g(v)-I'_n(v)v\qquad \mbox{and}\qquad \bar{g}(v)=g(v)-\bar{g}'(v)v,
\ee
after some algebra, one can write (\ref{middle_eqn}) as
\beq \label{middle_eqn2}
\lefteqn{[n\bar{g}'(v)-I'_n(v)]\frac{h(v)v}{H(v)}}\nonumber \\
    & = &
   \left[n+2-\frac{h'(v)H(v)}{h^2(v)}\right][\bar{g}'(v)-I'_n(v)]
+[(n+1)I'_n(v)-2\bar{g}'(v)]\frac{H(v)}{h(v)v}.
\eeq
Let us make the change of variables
\be \label{w_var}
w(v)=\frac{h(v)}{H(v)}v.
\ee
Note that $w(v)\equiv 1$ corresponds to the exponential c.d.f. (\ref{exp_type}). For simplicity, we write $w$ for $w(v)$.
After  multiplying both sides of (\ref{middle_eqn2}) by $h(v)v/H(v)$, using (\ref{w_var}) and
\[
\frac{h'(v)v}{h(v)}=\frac{w'}{w}v+w-1,
\]
it is not difficult to see that
(\ref{middle_eqn2}) is equivalent to

\be \label{main}
[n\bar{g}'(v)-I'_n(v)](w^2-1)=
\left[(n+1)(w-1)-\frac{w'}{w}v\right][\bar{g}'(v)-I'_n(v)].
\ee
Clearly, $w(v)\equiv 1$ is one solution, which corresponds to the exponential distribution. It remains to prove that this is the only solution.
Suppose $w(v)$ is a solution and there exists a value $v_1>0$ such that  $w(v_1)\ne 1$. We
want to reach a contradiction. Since $F$ is twice differentiable, we
have that $w(v)$ is continuous with respect to $v$ and hence $w(v)\ne 1$ for $v$ in an open
interval around $v_1$. Let
$v_0=\inf\{v|w(v)\ne 1\}$.
 Since $w(0+)=1$, we have
$v_0\ge l_F$.
We shall prove that $v_0=l_F$. Assume on contrary that $v_0>0$.
Then $w(v)=1$ if $0<v\le v_0$ and integration of (\ref{w_var}) implies that $h(v)$ is constant-valued in this interval. This contradicts the assumption (ii). Therefore $v_0=0$ and hence equation  (\ref{main}) holds for all $v>0$. Dividing (\ref{main}) by $w-1\ne 0$, we obtain
\beq \label{main2}
(w+1)[n\bar{g}'(v)-I'_n(v)]=
\left[n+1-\frac{w'v}{w(w-1)}\right][\bar{g}'(v)-I'_n(v)].
\eeq
It can be seen that
\[
\lim_{v\to 0+}\frac{w'v}{w(w-1)}=1.
\]
Also, using the L'Hopital's rule, one can find that
\[
\lim_{v\to 0+}\bar{g}'(v)=\frac{g'(0+)}{2}\qquad \mbox{and}\qquad   \lim_{v\to 0+}I'_n(v)=\frac{g'(0+)}{n+1}.
\]
Therefore, passing to the limit in (\ref{main2}) as $v\to 0+$, we obtain
\[
2\left(\frac{n}{2}-\frac{1}{n+1}\right)g'(0+)=n\left(\frac{1}{2}-\frac{1}{n+1}\right)g'(0+),
\]
which is not possible for $n\ge 2$.
This contradiction completes the proof.

\underline{Necessity}. This follows from the lemma setting $r=2$.

\end{document}